# First Integrals vs Limit Cycles

**Andrés G. García**

*Grupo de Investigación en Multifísica Aplicada (GIMAP), Universidad Tecnológica Nacional, 11 de Abril 461, Bahía Blanca, Buenos Aires, Argentina*

**Abstract**

This paper applies a recent result determining periodic orbits on the basis of first integrals, for Liénard systems. By solving a first order ODE with singularities, a crucial result is proved to locate intervals of single and isolated maximum amplitude's periodic orbits (limit cycles). With this result an upper bound for the number of limit cycles is provided. Some examples are presented along with conclusions and future work

1. **Introduction**

Physical systems' modeling ends up with nonlinear ordinary differential equations (ODE's) in the majority of the cases (see for instance [1]). However, an interesting and rich class of nonlinear systems is the case of nonlinear oscillators (see for instance [2] and [3]).

As it is well known, nonlinear ODE's cannot be integrated (solved) in closed-form except for some particular and special structures (this is the case of the simple pendulum [4]).

Moreover, many important properties can be studied nonlinear ODE's, however an important and common study regards the location and number of limit cycles (see for instance [5]).

In this case, few tools can be deployed to obtain conclusions ([1]–[7]). The problem of being isolated orbits that must be identified from the rest with the unique property of being periodic is neither trivial nor systematic.

At present, several methods to approximate (analyze) periodic orbits can be mentioned (see for instance [5] and the references therein):

- Homotopy perturbation

- Harmonic Balance

- Adomain Decomposition

- Variational Formulation

- Variational Iteration

- Pseudospectral method

- Rayleigh-Balance

- Energy-Balance

- Max-Min approach

- Amplitude-Frequency formulation

- Homotopy analysis

- Optimal Homotopy asymptotic method

- Pseudo-First Integrals (see [8])

In this paper, applying the first integral theorem recently proved in [7] to Liénard systems, an important location result is proved along with an upper bound for the number of limit cycles.

This paper is organized as follows: Section 2 presents the necessary machinery and preliminary results, Section 3 integrates an equivalent ODE for periodic orbits, Section 4 provides two important results: a section/counting sectors of limit cycles' existence and an upper bound for the number of limit cycles in Liénard systems, Section 5 presents some examples with known number of limit cycles and finally Section 6 presents some conclusions and future work.

**2. Preliminary results: Nonlinear oscillators' first integral**

A general nonlinear oscillator can be considered to be:

$$\ddot{x}(t) = f(x(t), \dot{x}(t)) \tag{1}$$

Where $f(x, \dot{x}): \mathbb{R}^2 \to \mathbb{R}$ is a nonlinear $C^0(\Re)$ function.

An important useful concept to distinguish a periodic orbit from other orbits has been introduced in [7], however a revision of first integral's definition is instructive:

**Definition (First Integral):** Given a system of first order ODE system (1), a first integral of the system is a continuously-differentiable function: $\Psi\ \Re \times \Re^n \to \Re^n$, which is not trivially constant but locally constant for any $x(t)$ solving (1):

$$\frac{d\Psi(t, x(t))}{dt} = 0 \tag{2}$$

From (2), a specific first integral for periodic orbits was defined and obtained in [7]

(the proof of this result is reproduced in this paper for the sake of completeness):

**Theorem 1:** A second order ODE (1) possess a periodic orbit: $\{x(0) = A \in \mathbb{R}^+, x(T) = A, \dot{x}(0) = 0\}$, if and only if there exists a function $\phi(x) \in C^1(\mathfrak{R})$, such that:

$$\begin{cases} \int_x^A f(x(t), \phi(x)) \cdot dx \leq 0 \\ \phi(A) = 0 \\ x \leq A \\ \dfrac{d\phi(x)}{dx} = \dfrac{f(x, \phi(x))}{\phi(x)} \end{cases} \quad (3)$$

*The proof is in the Appendix.*

### 3. First integral's applied to Liénard systems

An interesting analysis can be carried out for Liénard systems ([9]):

$$\begin{cases} \dot{x}(t) = y(t) - F(x(t)) \\ \dot{y}(t) = -x(t) \end{cases} \quad (4)$$

Where: $F(x) = a_n \cdot x^n + a_{n-1} \cdot x^{n-1} + \cdots + a_1 \cdot x$, for any $n \in \mathbb{N}$. This system can be readily transformed into a second order form:

$$\ddot{x}(t) = -x(t) - \frac{dF(x(t))}{dx(t)} \cdot \dot{x}(t) \quad (5)$$

Applying Theorem 1 (equation (3)) to equation (5), periodic orbits can be studied via the existence's solution of the following ODE:

$$\begin{cases} \dfrac{d\phi(x)}{dx} \cdot \phi(x) = -x - \dfrac{dF(x)}{dx} \cdot \phi(x) \\ \phi(A) = 0 \end{cases} \quad (6)$$

To divide the analysis, the following lemma is useful:

**Lemma 1:** Equation (6) possess continues solutions at least inside the intervals:

$$x \in [\bar{x}_{i-1}, \bar{x}_i], \quad i = 1, 2, \cdots n$$
$$\left.\frac{dF(x)}{dx}\right|_{\bar{x}_i} = 0$$

Where $\overline{x_n} = \infty$.

Note: This lemma is about solutions' continuity of (6) interior to the mentioned intervals, however nothing can be said about the borders: $\bar{x}_i$ , $i = 1, 2, \cdots n$.

*The proof is in the Appendix.*

**4 Maximal number of limit cycles in Liénard systems**

Lemma 1 provides the implicit suggestion to divide the analysis into intervals. However, periodic orbit's existence is equivalent to the existence of solutions crossing zero in equation (6):

$$Existence\ of\ a\ Periodic\ orbit\ (Amplitude\ A) \Leftrightarrow \begin{cases} \dfrac{d\phi(x)}{dx} = \dfrac{-x - \dfrac{dF(x)}{dx}\cdot\phi(x)}{\phi(x)} \\ \phi(A) = 0 \end{cases} \quad (7)$$

In other words, periodic orbits are only possible if an initial condition $\phi(\bar{x}_{i-1})$ exists to ensure a solution to equation (7) crossing zero. Moreover:

$$\lim_{x \to A} \frac{d\phi(x)}{dx} = \infty$$

Then the following lemma provides the maximum amplitudes limit cycles' location:

**Lemma 2 (Location of Limit cycles' maximum amplitudes):** Given a Liénard system (4) with a polynomial *F(x)* of degree *n*, an isolated periodic orbit (limit cycle) with maximal amplitude *A* (maximum amplitude in *x* coordinate), each interval:

$$x \in [\bar{x}_{i-1}, \bar{x}_i] \subseteq \mathbb{R}^+, \quad i = 1, 2, \cdots n$$

$$\left.\frac{dF(x)}{dx}\right|_{\bar{x}_i} = 0$$

contains at most one maximal $A$.

*The proof is in the Appendix.*

The immediate application of the lemma is the following upper bound for the number of limit cycles:

**Theorem 2 (Number of limit cycles):** Given a Liénard system (4) with a polynomial $F(x)$ of degree $n$, then the maximal number of limit cycles is $2 \cdot n - 1$.

*Proof:*

Applying Lemma 1, it straightforward to conclude that at most a single periodic orbit's amplitude $A$ can be counted inside each interval $[\bar{x}_{i-1}, \bar{x}_i]$, so at most n-1 intervals (amplitudes or periodic orbits) leaving the exception of possible amplitudes located at the borders: $\bar{x}_i, i = 1, \cdots n$ adding this number to the upper bound $n-1$, the maximal number of limit cycles can be $2 \cdot n - 1$. This completes the proof.

**4.1 Analysis of centers**

Systems (4) with centers rather than limit cycles must also be concluded from previous results.

In particular, the center's case is about non-isolated periodic orbits, so the conclusions in the proof of Lemma 2 precluding Bellman's possibility cannot be applied in the view of (A.6) (see Appendix):

$$\phi(x) = \pm \sqrt{\frac{(\bar{x}_{i-1}{}^2 - x^2 + \phi(\bar{x}_{i-1})^2)}{2} - \sum_{l=0}^{n-1} P_l(x) \cdot \frac{d^{n-l-1}y(x)}{dx^{n-l-1}}} \quad \cdots$$

If the origin (x=0) is a center, then at least for a dense interior set: $x \in \Omega \subseteq [\bar{x}_{i-1}, \bar{x}_i], \; i = 1, 2, \cdots n$:

$$\phi(x) = 0 = \pm \sqrt{\frac{(\bar{x}_{i-1}^2 - x^2 + \phi(\bar{x}_{i-1})^2)}{2} - \sum_{l=0}^{n-1} P_l(x) \cdot \frac{d^{n-l-1}y(x)}{dx^{n-l-1}}}$$

The conclusion is immediate, for a given initial condition $\phi(\bar{x}_{i-1})$, once reaching $\phi(A) = 0$, the solution to (6) remains null:

$$\phi(x) \sim 0 \quad (x \to A)$$

Confirming the ultimate behavior estimated by Bellman's in [11], pp. 101.

## 5. Examples

### Example 1

In order to apply Lemma 2 and Theorem 2, the well-known example in [9], pp. 345 is considered:

$$\begin{cases} \dot{x}(t) = y(t) - F(x(t)) \\ \dot{y}(t) = -x(t) \end{cases}$$

Where: $F(x) = a_3 \cdot x^3 + a_2 \cdot x^2 + a_1 \cdot x$. Then the following facts (see [9]) are known:

- If $a_1 \cdot a_3 > 0$, the system possess no closed orbits

- If $a_1 \cdot a_3 < 0$, the system possess a unique closed orbit

- If $a_1 = 0, \; a_3 \neq 0$, the system possess a weak attractor at the origin

- If $a_3 = 0, \; a_1 \neq 0$, the system possess an hy7perbolic attractor at the origin

- If $a_1 = 0$, $a_3 = 0$, the system possess a center at the origin

Lemma 2 provides regions of possible location of limit cycles:

$$\frac{dF(x)}{dx} = 3 \cdot a_3 \cdot x^2 + 2 \cdot a_2 \cdot x + a_1 = 3 \cdot a_3 \cdot \left(x + \frac{a_2}{3 \cdot a_3}\right)^2 + \left(\frac{3 \cdot a_1 \cdot a_3 - a_2{}^2}{3 \cdot a_3}\right)$$

- If $a_1 \cdot a_3 > 0$, then $\frac{dF(x)}{dx}$ possess no real roots, so at most a center at the origin is expected (no limit cycles)

- If $a_1 = 0$, $a_3 \neq 0$, then at most three limit cycles can be expected if: $\frac{a_2}{a_3} < 0$.

- If $a_3 = 0$, $a_1 \neq 0$, then at most three limit cycles can be expected if: $\frac{a_1}{a_2} < 0$.

- If $a_1 = 0$, $a_3 = 0$, then at most one limit cycle can be expected.

As for Theorem 2, the upper bound on the number of limit cycles indicates: $2 \cdot 3 - 1 = 5$.

**Example 2**

In [16], several Liénard systems were considered, in particular:

$$F(x) = x^5 - \mu \cdot x^3 + x$$

For this system is known (see [16]):

- If $\mu > 2.5$, then the system has exactly two limit cycles

- If $\mu < 2$, then the system possess no limit cycles

- If $\mu \in [2, 2.5]$, then the system exhibits a bifurcation zone

Lemma 2 leads:

$$\frac{dF(x)}{dx} = 5 \cdot x^4 - \mu \cdot 3 \cdot x^2 + 1 \Rightarrow x^2 = \frac{3 \cdot \mu \pm \sqrt{9 \cdot \mu^2 - 20}}{10}$$

The condition for existence of positive roots yields: $\mu > \frac{2 \cdot \sqrt{5}}{3} \cong 1.4907$. This condition agree with the known value $\mu > 2$, since Lemma 2 is a necessary condition. Notice that the ultimate case $\mu > \frac{2 \cdot \sqrt{5}}{3}$ will lead at most 5 limit cycles.

Also in [16] another system with known number of limit cycles was studied:

$$F(x) = 0.32 \cdot x^5 - \frac{4}{3} \cdot x^3 + 0.8 \cdot x$$

In this case, Lemma 2 provides at most 5 limit cycles (two positive roots of $\frac{dF(x)}{dx}$). In this case is instructive to visualize the location of the maximum amplitude of the two limit cycles to verify that they are in fact inside the intervals prescribed by Lemma 2 (see Figure 1):

$$\frac{dF(x)}{dx} = 0.32 \cdot 5 \cdot x^4 - 4 \cdot x^2 + 0.8 = 0 \Rightarrow x = \{1.5102, 0.4682\}$$

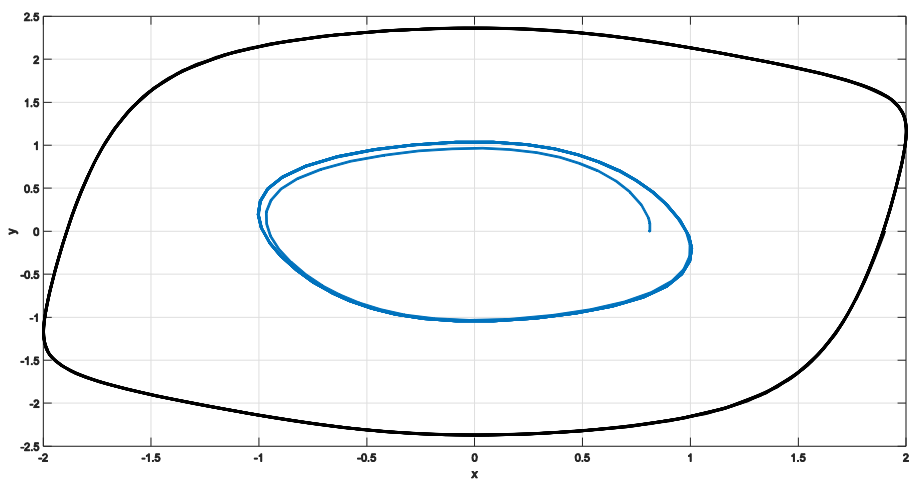

Figure 1: Numerical analysis of limit cycles

**Example 3**

The last example considers the system studied in [17]:

$$F(x) = (x^2 - 1)^2 \cdot (c \cdot e \cdot x + 1) \cdot \left(x^2 + e \cdot x + \frac{1}{8}\right)$$

With $\{c, e\}\mathbb{R}^+$. Then Lemma 2 leads:

$$\frac{dF(x)}{dx} = (x^2 - 1)^2 \cdot (c \cdot e \cdot x + 1) \cdot \left(x^2 + e \cdot x + \frac{1}{8}\right) = 0$$

As analyzed in [17] for small values of *e* the critical points (derivatives of *F(x)* equal to zero) are located at: $\left\{r, \frac{1}{2} + O(e), 1\right\}$, where $0 < r < \frac{1}{2}$. In this case, Lemma 2 allows at most 7 limit cycles.

## 6. Conclusions

In this paper we suggest a novel way to count and locate interval of (possible) existence of isolated periodic orbits (limit cycles) in Liénard systems. The main theorem is a consequence of the application of the first integral concept for periodic orbits presented in [7].

Integrating an equivalent first order ODE with singularities, a lemma enclosing maximal amplitudes of limit cycles allowed the derivation of an upper bound for limit cycles.

While several examples were presented resorting to Liénard systems with known number of limit cycles, numerical integration also showed the application of the main results.

As a future work, numerical and hardware (analog and digital) implementation of method based on piece-wise linear ideas will extract more information about the shape, period and exact amplitude of limit cycles in Liénard as well as general nonlinear oscillators.

**Acknowledgments:** The work is supported by Universidad Tecnológica Nacional.

## 7. APPENDIX

**<u>Proof of Theorem 1:</u>**

Necessity:

If there exists $\phi(x) \in C^0(\Re)$, then:

$$\ddot{x}(t) = f(x(t), \phi(x(t)))$$

According to [12], pp. 26-27:

$$\frac{d\dot{x}(t)}{dx} = -\frac{f(x(t), \phi(x(t)))}{\dot{x}(t)} \Leftrightarrow d\dot{x}^2(t) = -f(x(t), \phi(x(t))) \cdot dx$$

Equivalently:

$$\frac{\dot{x}^2(t)}{2} = -\int_A^x f(x(t), \phi(x(t))) \cdot dx$$

Existence of solutions are possible only if: $\dfrac{\dot{x}^2(t)}{2} = -\int_A^x f(x(t), \phi(x(t))) \cdot dx \geq 0$, which

is ascertained by hypothesis. On the other hand and taking into account that

$\ddot{x}(t) = f(x(t), \phi(x(t)))$ does not depend on $\dot{x}(t)$, the reasoning on [13] is applied thus

proving the existence of periodic orbits: $\dot{x}(t) = \phi(x(0) = A) = 0$.

Sufficiency:

If periodic orbits with period $T$ do exists, then: $\ddot{x}(t) = f(x, \dot{x})$, $\dot{x}(0) = \dot{x}(T) = 0$.

Equivalently:

$$\int_0^T f(x(t), \dot{x}(t)) \cdot dt = 0$$

By contradiction, let's assume that do not exists $\phi(x) \in C^0(\Re)$ such that $\dot{x}(t)$ and

$x(t)$ are related, then an asymptotic expansion for an arbitrary bounded function

$\varepsilon(t) \in C1(\Re)$ (see for instance [14]) can be carried out:

$$f(x,\dot{x}) \approx f\big(x(t)+\eta \cdot \varepsilon(t), \dot{x}(t)+\eta \cdot \dot\varepsilon(t)\big) + \frac{\partial f(x,\dot{x})}{\partial x}\bigg|_{\substack{x(t)+\eta\cdot\varepsilon(t)\\ \dot{x}(t)+\eta\cdot\dot\varepsilon(t)}} \cdot \big(x-(x+\eta\cdot\varepsilon(t))\big)$$

$$+ \frac{\partial f(x,\dot{x})}{\partial \dot{x}}\bigg|_{\substack{x(t)+\eta\cdot\varepsilon(t)\\ \dot{x}(t)+\eta\cdot\dot\varepsilon(t)}} \cdot \big(\dot{x}-(\dot{x}+\eta\cdot\dot\varepsilon(t))\big) \qquad (\eta \to 0)\,\forall\, \varepsilon(t)\in C^1(\Re)$$

Integrating:

$$\int_0^T f(x,\dot{x})\cdot dt = 0 \approx \int_0^T f\big(x(t)+\eta\cdot\varepsilon(t), \dot{x}(t)+\eta\cdot\dot\varepsilon(t)\big)\cdot dt +$$

$$-\eta \cdot \int_0^T \left[ \frac{\partial f(x,\dot{x})}{\partial x}\bigg|_{\substack{x(t)+\eta\cdot\varepsilon(t)\\ \dot{x}(t)+\eta\cdot\dot\varepsilon(t)}} \cdot \varepsilon(t) + \frac{\partial f(x,\dot{x})}{\partial \dot{x}}\bigg|_{\substack{x(t)+\eta\cdot\varepsilon(t)\\ \dot{x}(t)+\eta\cdot\dot\varepsilon(t)}} \cdot \dot\varepsilon(t) \right]\cdot dt \qquad (\eta \to 0)\,\forall\, \varepsilon(t)\in C^1(\Re)$$

Then:

$$\int_0^T f\big(x(t)+\eta\cdot\varepsilon(t), \dot{x}(t)+\eta\cdot\dot\varepsilon(t)\big)\cdot dt \approx -\eta \cdot \int_0^T \left[ \frac{\partial f(x,\dot{x})}{\partial x}\bigg|_{\substack{x(t)+\eta\cdot\varepsilon(t)\\ \dot{x}(t)+\eta\cdot\dot\varepsilon(t)}} \cdot \varepsilon(t) + \frac{\partial f(x,\dot{x})}{\partial \dot{x}}\bigg|_{\substack{x(t)+\eta\cdot\varepsilon(t)\\ \dot{x}(t)+\eta\cdot\dot\varepsilon(t)}} \cdot \dot\varepsilon(t) \right]\cdot dt \qquad (\eta \to 0)$$

For one side:

$$\int_0^T f\big(x(t)+\eta\cdot\varepsilon(t), \dot{x}(t)+\eta\cdot\dot\varepsilon(t)\big)\cdot dt \approx \int_0^T f\big(x(t),\dot{x}(t)\big)\cdot dt = 0 \qquad (\eta \to 0)$$

Then:

$$-\eta \cdot \int_0^T \left[ \frac{\partial f(x,\dot{x})}{\partial x}\bigg|_{\substack{x(t)+\eta\cdot\varepsilon(t)\\ \dot{x}(t)+\eta\cdot\dot\varepsilon(t)}} \cdot \varepsilon(t) + \frac{\partial f(x,\dot{x})}{\partial \dot{x}}\bigg|_{\substack{x(t)+\eta\cdot\varepsilon(t)\\ \dot{x}(t)+\eta\cdot\dot\varepsilon(t)}} \cdot \dot\varepsilon(t) \right]\cdot dt \approx 0 \qquad (\eta \to 0)\,\forall\, \varepsilon(t)\in C^1(\Re)$$

Since these asymptotic equivalence is valid for any arbitrary function $\varepsilon(t) \in C1(\Re)$, if

$$\int_0^T \left[ \left.\frac{\partial f(x,\dot{x})}{\partial x}\right|_{\substack{x(t)+\eta\cdot\varepsilon(t) \\ \dot{x}(t)+\eta\cdot\varepsilon(t)}} \cdot \varepsilon(t) + \left.\frac{\partial f(x,\dot{x})}{\partial \dot{x}}\right|_{\substack{x(t)+\eta\cdot\varepsilon(t) \\ \dot{x}(t)+\eta\cdot\varepsilon(t)}} \cdot \dot{\varepsilon}(t) \right] \cdot dt \quad \text{is not asymptotic equivalent to}$$

zero, there should exists a maxima or tend to infinity:

<u>T<∞</u>:

$$\int_0^T \left[ \left.\frac{\partial f(x,\dot{x})}{\partial x}\right|_{\substack{x(t)+\eta\cdot\varepsilon(t) \\ \dot{x}(t)+\eta\cdot\varepsilon(t)}} \cdot \varepsilon(t) + \left.\frac{\partial f(x,\dot{x})}{\partial \dot{x}}\right|_{\substack{x(t)+\eta\cdot\varepsilon(t) \\ \dot{x}(t)+\eta\cdot\varepsilon(t)}} \cdot \dot{\varepsilon}(t) \right] \cdot dt \leq$$

$$\leq \int_0^T \left[ \left.\frac{\partial f(x,\dot{x})}{\partial x}\right|_{\substack{x(t)+\eta\cdot\varepsilon(t) \\ \dot{x}(t)+\eta\cdot\varepsilon(t)}} \cdot \bar{\varepsilon}(t) + \left.\frac{\partial f(x,\dot{x})}{\partial \dot{x}}\right|_{\substack{x(t)+\eta\cdot\varepsilon(t) \\ \dot{x}(t)+\eta\cdot\varepsilon(t)}} \cdot \dot{\bar{\varepsilon}}(t) \right] \cdot dt, \quad \forall \varepsilon(t) \in C^1(\Re)$$

Where $\bar{\varepsilon}(t)$ is the particular function ε(t) that produces the maximum. This bound is valid for arbitrary function ε(t)∈C1(ℜ), so it is also valid for

$\varepsilon(t) = \rho \cdot \bar{\varepsilon}(t), \quad \forall \rho \in \Re$, then is valid for $\rho > 1$, so the bound is not the maxima.

This conclusion explains that:

$$\int_0^T \left[ \left.\frac{\partial f(x,\dot{x})}{\partial x}\right|_{\substack{x(t)+\eta\cdot\varepsilon(t) \\ \dot{x}(t)+\eta\cdot\varepsilon(t)}} \cdot \varepsilon(t) + \left.\frac{\partial f(x,\dot{x})}{\partial \dot{x}}\right|_{\substack{x(t)+\eta\cdot\varepsilon(t) \\ \dot{x}(t)+\eta\cdot\varepsilon(t)}} \cdot \dot{\varepsilon}(t) \right] \cdot dt = 0 \quad \forall \varepsilon(t) \in C^1(\Re)$$

Integrating by parts and because of the fundamental lemma of calculus of variations (see [15], pp.9 Lemma 1):

$$\frac{d}{dt}\left(\frac{\partial f(x,\dot{x})}{\partial \dot{x}}\right) - \frac{\partial f(x,\dot{x})}{\partial x} = 0 \Leftrightarrow \varphi(x,\dot{x}) = 0$$

Valid only if $\dot{x}(t)$ and $x(t)$ are not a function each other, in other words, it is valid only if does not exist $\phi(x) \in C^0(\Re)$, since derivatives are taken independently for $\{x(t), \dot{x}(t)\}$.

Finally, the implicit function theorem (see for instance [6], pp. 50) implies:

$$\frac{d}{dt}\left(\frac{\partial f(x,\dot{x})}{\partial \dot{x}}\right) - \frac{\partial f(x,\dot{x})}{\partial x} = 0 \Leftrightarrow \varphi(x,\dot{x}) = 0 \Leftrightarrow \frac{\partial \varphi(x,\dot{x})}{\partial x} = 0$$

For every periodic trajectory of $\ddot{x}(t) = f(x,\dot{x})$, $f : \Re^2 \to \Re$. Then:

$$\varphi(x,\dot{x}) = 0 \Leftrightarrow \frac{\partial \varphi(x,\dot{x})}{\partial x} = 0 \Leftrightarrow \dot{\varphi}(x,\dot{x}) = \frac{\partial \varphi(x,\dot{x})}{\partial x} \cdot \dot{x} + \frac{\partial \varphi(x,\dot{x})}{\partial \dot{x}} \cdot \ddot{x} = 0$$

This is:

$$\dot{\varphi}(x,\dot{x}) = \underbrace{\frac{\partial \varphi(x,\dot{x})}{\partial x}}_{0} \cdot \dot{x} + \frac{\partial \varphi(x,\dot{x})}{\partial \dot{x}} \cdot \ddot{x} = 0 \Leftrightarrow \frac{\partial \varphi(x,\dot{x})}{\partial \dot{x}} \cdot f(x,\dot{x}) = 0$$

In summary, every periodic trajectory of $\ddot{x}(t) = f(x,\dot{x})$, $f : \Re^2 \to \Re$, satisfies:

$$\begin{cases} \dfrac{\partial \varphi(x,\dot{x})}{\partial x} = 0 \\ \dfrac{\partial \varphi(x,\dot{x})}{\partial \dot{x}} = 0 \end{cases}$$

A phase-portrait can be depicted in Figure A.1:

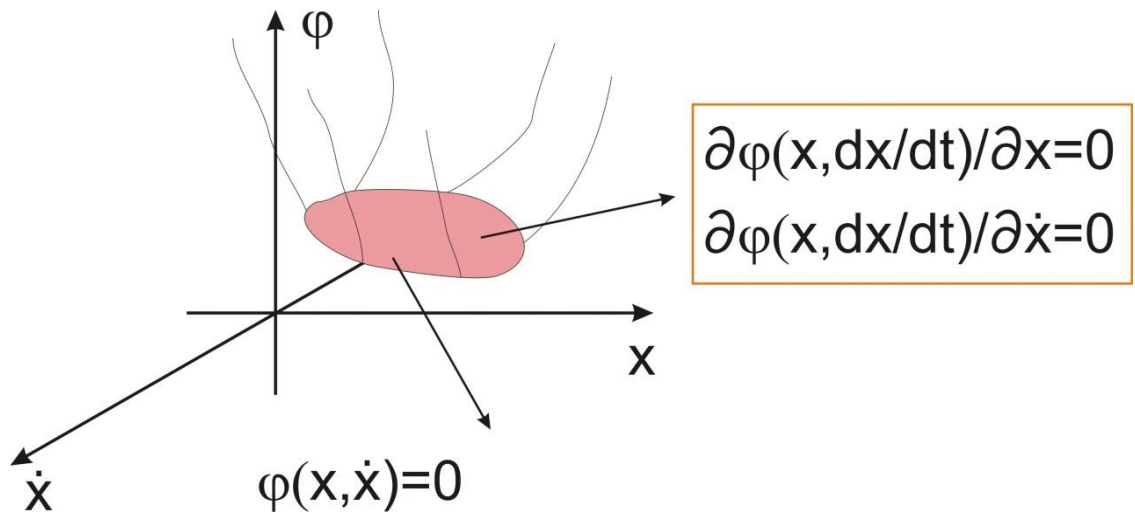

Figure A.1: 3-D Phase-Portrait assuming the non-existence of $\phi(x)$

This is equivalent to a phase-portrait indicated in Figure A.2:

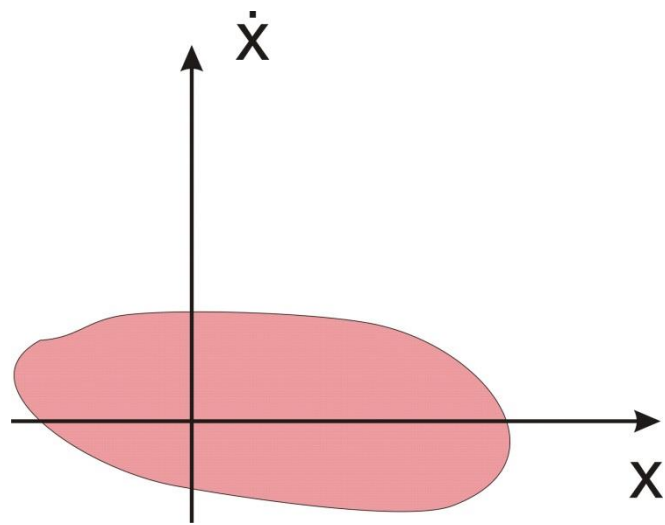

Figure A.2: 2-D Phase-Portrait assuming the non-existence of $\phi(x)$

This is a contradiction to the non-existence hypothesis for $\phi(x)$. Once this existence is proved for periodic orbits, a similar reasoning as in the previous (necessity) part can be carried out. This completes the proof.

T→∞:

$$-\eta \cdot \int_0^\infty \left[ \frac{d}{dt}\left( \frac{\partial f(x,\dot{x})}{\partial \dot{x}} \right) - \frac{\partial f(x,\dot{x})}{\partial x} \right] \cdot \varepsilon(t) \cdot dt \approx 0 \quad (\eta \to 0) \; \forall \; \varepsilon(t) \in C^1(\Re)$$

In case that the integral possess a maximum, the reasoning is as in previous case. So, the only missing analysis is: $\int_0^\infty \left[ \frac{d}{dt}\left( \frac{\partial f(x,\dot{x})}{\partial \dot{x}} \right) - \frac{\partial f(x,\dot{x})}{\partial x} \right] \cdot \varepsilon(t) \cdot dt \to \infty$. This

possibility may occur for some bounded function $\varepsilon(t) \in C^1(\Re)$:

$$\int_0^\infty \frac{d}{dt}\left( \frac{\partial f(x,\dot{x})}{\partial \dot{x}} \right) \cdot \varepsilon(t) \cdot dt - \int_0^\infty \frac{\partial f(x,\dot{x})}{\partial x} \cdot \varepsilon(t) \cdot dt \to \infty$$

Two cases must be considered:

- $\int_0^\infty \frac{d}{dt}\left( \frac{\partial f(x,\dot{x})}{\partial \dot{x}} \right) \cdot \varepsilon(t) \cdot dt \to \infty$

- $\int_0^\infty \frac{\partial f(x,\dot{x})}{\partial x} \cdot \varepsilon(t) \cdot dt \to \infty$

The first case is not possible taking into account the hypothesis of periodic orbits $x(t)$, if $\varepsilon(t)$ is a function of $x(t)$, so $\varepsilon(t)$ is periodic with the integral of a periodic orbit in a period is indeed zero. Otherwise, if $\varepsilon(t)$ is not a function of $x(t)$:

$$\varepsilon(t) \cdot \left[ \left. \frac{\partial f(x,\dot{x})}{\partial \dot{x}} \right|_{x(0)} - \left. \frac{\partial f(x,\dot{x})}{\partial \dot{x}} \right|_{x(\infty)} \right] = 0$$

The second case requires a more detailed analysis, integrating by parts:

$$\int_0^\infty \frac{\partial f(x,\dot{x})}{\partial x} \cdot \varepsilon(t) \cdot dt = \varepsilon(t) \cdot \underbrace{\int_0^\infty \frac{\partial f(x,\dot{x})}{\partial x} \cdot dt}_{Periodic} - \int_{\varepsilon(0)}^{\varepsilon(\infty)} \int_0^t \frac{\partial f(x,\dot{x})}{\partial x} \cdot d\sigma \cdot d\varepsilon(t)$$

Then:

$$\int\limits_{0}^{\infty} \frac{\partial f(x,\dot{x})}{\partial x} \cdot \varepsilon(t) \cdot dt = -\int\limits_{\varepsilon(0)}^{\varepsilon(\infty)} \int\limits_{0}^{t} \frac{\partial f(x,\dot{x})}{\partial x} \cdot d\sigma \cdot d\varepsilon(t) \to \infty$$

As in the previous case, the valid possibility is $\varepsilon(t)$ not a function of $x(t)$:

$$\int\limits_{0}^{\infty} \frac{\partial f(x,\dot{x})}{\partial x} \cdot \varepsilon(t) \cdot dt = -\int\limits_{0}^{t} \frac{\partial f(x,\dot{x})}{\partial x} \cdot d\sigma \cdot \underbrace{\int\limits_{\varepsilon(0)}^{\varepsilon(\infty)} d\varepsilon(t)}_{Bounded} \to \infty$$

Since $\varepsilon(t)$ is bounded, this possibility is not allowed. Once the finiteness of the integral was proven, the reasoning is as in previous case. This completes the proof.

**Proof of Lemma 1:**

Integrating (6) from $x_0$ to $x$:

$$-\frac{\phi(x)^2}{2} + \frac{\phi(x_0)^2}{2} = \frac{(x^2 - x_0^2)}{2} + \int_{x_0}^{x} \frac{dF(\sigma)}{d\sigma} \cdot \phi(\sigma) \cdot d\sigma \tag{A.1}$$

To solve the integral $\int_{x_0}^{x} \frac{dF(\sigma)}{d\sigma} \cdot \phi(\sigma) \cdot d\sigma$ in (A.1), the polynomial structure of $F(x)$ is exploited:

$$\int_{x_0}^{x} \frac{dF(\sigma)}{d\sigma} \cdot \phi(\sigma) \cdot d\sigma = \int_{x_0}^{x} \sum_{i=1}^{n} a_i \cdot i \cdot \sigma^{i-1} \cdot \phi(\sigma) \cdot d\sigma$$

Integrating by parts each term:

$$\int_{x_0}^{x} \sigma^{i-1} \cdot \phi(\sigma) \cdot d\sigma$$

$$= \sigma^{i-1} \cdot \int \phi(\sigma_1) \cdot d\sigma_1 \bigg|_{x_0}^{x} - (i-1) \cdot \int_{x_0}^{x} \sigma^{i-2} \cdot \int \phi(\sigma_1) \cdot d\sigma_1 \cdot d\sigma$$

Defining:

$$Z_1(x) = \int_{x_0}^{x} \phi(\sigma_1) \cdot d\sigma_1$$

Then:

$$\int_{x_0}^{x} \sigma^{i-1} \cdot \phi(\sigma) \cdot d\sigma = -x^{i-1} \cdot Z_1(x) - (i-1) \cdot \int_{x_0}^{x} \sigma^{i-2} \cdot \int \phi(\sigma_1) \cdot d\sigma_1 \cdot d\sigma$$

Integrating again by parts:

$$\int_{x_0}^{x} \sigma^{i-1} \cdot \phi(\sigma) \cdot d\sigma$$

$$= -x^{i-1} \cdot Z_1(x) - (i-1)$$

$$\cdot \left[ x^{i-2} \cdot Z_2(x) - (i-2) \cdot \int_{x_0}^{x} \sigma^{i-3} \cdot \int \int \phi(\sigma_1) \cdot d\sigma_1 \cdot d\sigma_2 \cdot d\sigma \right]$$

Where:

$$Z_2(x) = \int_{x_0}^{x} \int \phi(\sigma_1) \cdot d\sigma_1 \cdot d\sigma$$

Following this procedure $i$ times:

$$\int_{x_0}^{x} \sigma^{i-1} \cdot \phi(\sigma) \cdot d\sigma$$

$$= x^{i-1} \cdot Z_1(x) - (i-1) \cdot x^{i-2} \cdot Z_2(x) + (i-2) \cdot x^{i-3} \cdot Z_3(x) - \cdots$$

$$+ x \cdot Z_{i-1}(x) - Z_i(x)$$

In a compact form:

$$\int_{x_0}^{x} \sigma^{i-1} \cdot \phi(\sigma) \cdot d\sigma = -\sum_{j=0}^{i} (-1)^j \cdot x^{i-j-1} \cdot \binom{i}{i-j} \cdot j! \cdot \frac{dZ_i(x)^{i-j-1}}{dx^{i-j-1}} \quad (A.2)$$

Where:

$$\frac{d^p Z_i(x)}{dx^p} = Z_{i-p}(x), i, p = 1, \cdots n$$
$$Z_0(x) = \phi(x) \quad (A.3)$$
$$Z_1(x) = \int_{x_0}^{x} \phi(\sigma_1) \cdot d\sigma_1$$

From (A.2) using (A.3):

$$-\frac{\phi(x)^2}{2} + \frac{\phi(x_0)^2}{2} = \frac{(x^2 - x_0^2)}{2} - \sum_{i=1}^{n} a_i \cdot i \cdot \sum_{j=0}^{i}(-1)^j \cdot x^{i-j-1} \cdot \binom{i}{i-j} \cdot j! \cdot \frac{dZ_i(x)^{i-j-1}}{dx^{i-j-1}}$$

To obtain a compact differential equation form, notice:

$$\begin{array}{cccc}
i - j - 1 = 0, & i - j = 1, & j = i - 1, & i \geq 1 \\
i - j - 1 = 1, & i - j = 2, & j = i - 2, & i \geq 2 \\
& \vdots & & \\
i - j - 1 = i - 1, & i - j = i, & j = 0, & i \geq i
\end{array}$$

Then:

$$-\frac{\phi(x)^2}{2} + \frac{\phi(x_0)^2}{2} = \frac{(x^2 - x_0^2)}{2} - \sum_{i=1}^{n} a_i \cdot i \cdot (-1)^{i-1} \cdot x^0 \cdot \binom{i}{1} \cdot (i-1)! \cdot$$

$$Z_i(x) + \sum_{i=2}^{n} a_i \cdot i \cdot (-1)^{i-2} \cdot x^1 \cdot \binom{i}{2} \cdot (i-2)! \cdot \frac{dZ_i(x)}{dx} + \sum_{i=3}^{n} a_i \cdot i \cdot (-1)^{i-3} \cdot x^2 \cdot$$

$$\binom{i}{3} \cdot (i-3)! \cdot \frac{d^2 Z_i(x)}{dx^2} + \cdots + a_n \cdot n \cdot (-1)^0 \cdot x^{n-1} \cdot \binom{i}{n} \cdot (0)! \cdot \frac{d^{n-1} Z_i(x)}{dx^{n-1}}$$

Considering (A.3):

$$-\frac{\phi(x)^2}{2} + \frac{\phi(x_0)^2}{2} = \frac{(x^2 - x_0^2)}{2} - \sum_{k=1}^{n} a_k \cdot k \cdot x^{k-1} \cdot Z_1(x) + \sum_{k=1}^{n-1} a_{k+1} \cdot (k+1) \cdot$$

$$(-1) \cdot (k-1) \cdot 1! \cdot x^{k-1} \cdot Z_2(x) + \sum_{k=1}^{n-2} a_{k+2} \cdot (k+2) \cdot (-1)^2 \cdot \binom{k+2}{k} \cdot 2! \cdot x^2 \cdot$$

$$Z_3(x) + \cdots + a_n \cdot n \cdot (-1)^{n-1} \cdot \binom{n}{1} \cdot (n-1)! \cdot Z_n(x)$$

(11)

Compactly:

$$-\frac{\phi(x)^2}{2} + \frac{\phi(x_0)^2}{2} = \frac{(x^2 - x_0^2)}{2} - \sum_{l=0}^{n-1} P_l(x) \cdot Z_{l+1}(x)$$

$$P_l(x) = \sum_{k=1}^{n-l} a_{k+l} \cdot (k+l) \cdot (-1)^l \cdot \binom{k+l}{k} \cdot l! \cdot x^{k-1}, \quad l = 0, \cdots, n-1$$

(A.4)

Defining:

$$Z_n(x) = y(x)$$
$$Z_{n-1}(x) = \frac{dy(x)}{dx}$$
$$Z_{n-1}(x) = \frac{d^2y(x)}{dx^2}$$
$$\vdots$$
$$Z_1(x) = \frac{d^{n-1}y(x)}{dx^{n-1}}$$

Equation (A.4) takes the form:

$$\sum_{l=0}^{n-1} P_l(x) \cdot \frac{d^{n-l-1}y(x)}{dx^{n-l-1}} = -\frac{\phi(x)^2}{2} + \frac{\phi(x_0)^2}{2} - \frac{(x^2 - x_0^2)}{2} \tag{A.5}$$

This ODE is in fact a Linear Time Varying ODE that can be rewritten in a vector form:

$$\frac{d\omega(x)}{dx} = \bar{A}(x) \cdot \omega + \begin{bmatrix} 0 \\ \vdots \\ 0 \\ 1 \end{bmatrix} \cdot \left( \frac{-x^2 + x_0^2 - \phi(x)^2 + \phi(x_0)^2}{2 \cdot P_0(x)} \right) \tag{13}$$

Where:

$$\bar{A}(x) = \begin{bmatrix} 0 & 1 & 0 & \cdots & 0 \\ 0 & 0 & 1 & \cdots & 0 \\ \vdots & & & \ddots & 1 \\ \frac{P_{n-1}(x)}{P_0(x)} & \frac{P_{n-2}(x)}{P_0(x)} & \cdots & & \frac{P_1(x)}{P_0(x)} \end{bmatrix}, \quad \omega = \begin{bmatrix} y(x) \\ \frac{dy(x)}{dx} \\ \frac{d^2y(x)}{dx^2} \\ \vdots \\ \frac{d^{n-2}y(x)}{dx^{n-2}} \end{bmatrix},$$

Existence and continuity of solutions can only assured if the ODE is at least $C^1(\mathbb{R}^{n-1})$:

$$\begin{matrix} x \in [\bar{x}_{i-1}, \bar{x}_i] \\ P_0(\bar{x}_i) = 0 \end{matrix}, i = 1, \cdots, n-1$$

Noticing that $\frac{dF(x)}{dx} = P_0(x)$. This completes the proof.

**Proof of Lemma 2:**

Let's suppose the existence of two periodic orbits (isolated) inside the interval $[\bar{x}_{i-1}, \bar{x}_i]$ for:

$$\begin{cases} \frac{d\phi(x)}{dx} \cdot \phi(x) = -x - \frac{dF(x)}{dx} \cdot \phi(x) \\ \phi(A) = 0 \end{cases} \quad (A.5)$$

From (A.4) specializing for $x_0 = \bar{x}_{i-1}$:

$$\phi(x) = \pm\sqrt{\frac{(\bar{x}_{i-1}^2 - x^2 + \phi(\bar{x}_{i-1})^2)}{2} + \sum_{l=0}^{n-1} P_l(x) \cdot \frac{d^{n-l-1}y(x)}{dx^{n-l-1}}} \quad (A.6)$$

Singularities may occur if and only if:

$$\frac{(\bar{x}_{i-1}^2 - x^2 + \phi(\bar{x}_{i-1})^2)}{2} + \sum_{l=0}^{n-1} P_l(x) \cdot \frac{d^{n-l-1}y(x)}{dx^{n-l-1}}, \forall x > A$$

However, this conditions implies the non-continuity of the solutions (A.6) for $x$ bigger than $A$. To preclude this possibility, continuity of (A.5) for $x<A$ is invoked, so the solutions are at least continuous for $x<A$ (see [2], pp.3-4), then if the solutions to (A.5) do not continue for $x>A$, the ultimate behavior of solutions should finish at $x=A$.

To this aim, a theorem in [11], pp. 101 is particularly useful:

*Any solution of* $\frac{d\phi(x)}{dx} = \frac{P(\phi,x)}{Q(\phi,x)}$, *continuous for* $x \geq \bar{x}_{i-1}$ *satisfies (ultimately):*

$$\phi(x) \sim a \cdot x^b \cdot e^{L(x)} 0 \text{ or } \phi(x) \sim a \cdot x^b \cdot (\log x)^{\frac{1}{c}}$$

*Where* $b \in \mathbb{R}$ *and* $c \in \mathbb{N}$ *are constants and* $\{P(\phi, x), Q(\phi, x), L(x)\}$ *polynomials.*

In the view of this result, only three possibilities arise:

$\underline{L(A) = -\infty}$

This condition cannot be satisfied for any finite $A$ considering the polynomial nature of $L(x)$.

$\underline{(\log A)^{\frac{1}{c}}}$

This condition implies $\log A = 0$, or in other words: $A=1$, leaving a particular case but also proving the single existence of periodic isolated orbits inside the intervals $[\bar{x}_{i-1}, \bar{x}_i]$.

$\underline{A^b = 0}$

This possibility implies:

$$\phi(x) \sim 0 \quad (x \to A)$$

By definition of asymptotic equivalence (see [10], pp. 10), this conclusion means:

$$\lim_{x \to A^-} \frac{\phi(x)}{x - A} = 1$$

According to (A.6):

$$\lim_{x \to A^-} \frac{\phi(x)}{x - A} = 1$$

$$\Leftrightarrow \lim_{x \to A^-} \frac{\pm \sqrt{\frac{(\bar{x}_{i-1}^2 - x^2 + \phi(\bar{x}_{i-1})^2)}{2} + \sum_{l=0}^{n-1} P_l(x) \cdot \frac{d^{n-l-1} y(x)}{dx^{n-l-1}}}}{x - A} = 1$$

Moreover, this limit is only possible if:

$$lim_{x \to A^-} \sqrt{\frac{(\bar{x}_{i-1}^2 - x^2 + \phi(\bar{x}_{i-1})^2)}{2} + \sum_{l=0}^{n-1} P_l(x) \cdot \frac{d^{n-l-1}y(x)}{dx^{n-l-1}}} = 0$$

Applying L'Hospital rule:

$$lim_{x \to A^-} \frac{\pm(-x)}{\sqrt{\frac{(\bar{x}_{i-1}^2 - x^2 + \phi(\bar{x}_{i-1})^2)}{2} + \sum_{l=0}^{n-1} P_l(A) \cdot \frac{d^{n-l-1}y(x)}{dx^{n-l-1}}\bigg|_{x=A}}} = \infty \neq 1$$

The obvious conclusion is the continuation of solutions for $x \geq A$. The possible case follows:

- Collision of solutions

    In this case, if a second periodic orbit found its amplitude inside the interval $[\bar{x}_{i-1}, \bar{x}_i]$, then Figure A.3 depicts the scenario. Clearly there exists at least a single point possessing multiple solutions for a given initial condition, this is a contradiction in the view of the ODE's continuity.

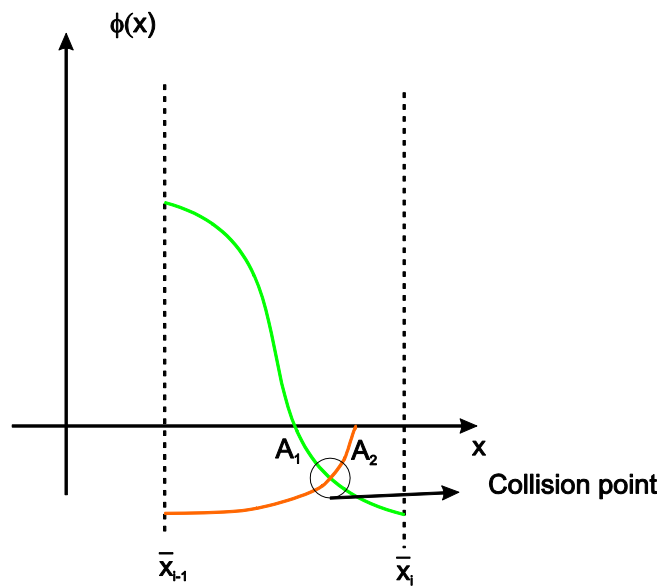

Figure A.3: The first possibility

- Rolle's theorem

In this case, the solution $\phi(x)$ is a continuous function crossing zero in two points (Figure A.4), so Rolle's theorem ensures the existence of an interior point $x^*$ such that:

$$\left.\frac{d\phi(x)}{dx}\right|_{x=x^*} = 0 \Leftrightarrow -x^* - \left.\frac{dF(x)}{dx}\right|_{x=x^*} \cdot \phi(x^*) = 0$$

Moreover, inside the interval $[\bar{x}_{i-1}, \bar{x}_i]$, according to Lemma 1: $P_0(x) \neq 0$ or equivalently $\frac{dF(x)}{dx} \neq 0$, so the sign of $\phi(x^*) = \frac{-x^*}{\left.\frac{dF(x)}{dx}\right|_{x=x^*}}$ remains unchanged, however the opposite sign in $\phi(x)$ also exists as a solution of (A.5) in the view of its periodic orbit' existence equivalence, this is a contradiction and completes the proof.

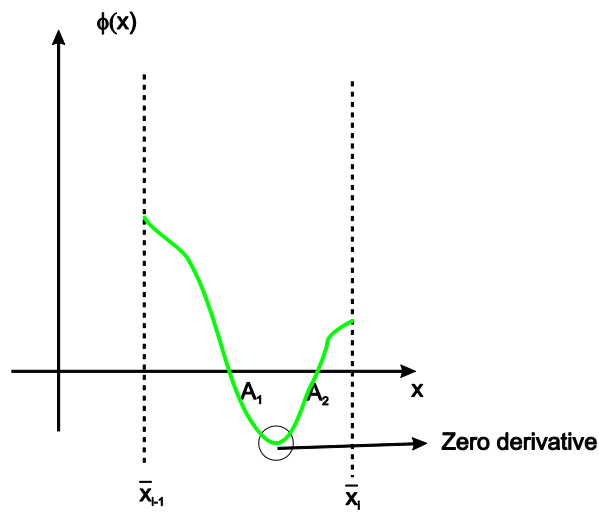

Figure A.4: The second possibility.